\documentclass[10pt]{amsart}
\usepackage[leqno]{amsmath}

\newtheorem{thm}{Theorem}
\newtheorem{cor}[thm]{Corollary}
\newtheorem{lem}[thm]{Lemma}

\theoremstyle{definition}
\newtheorem{definition}[thm]{Definition}
\newtheorem{remark}[thm]{Remark}
\numberwithin{thm}{section}
\numberwithin{equation}{section}

\newcommand{\EQ}[1]{\eqref{eq:#1}}
\newcommand{\LEM}[1]{Lemma~\ref{lem:#1}}
\newcommand{\DEF}[1]{Definition~\ref{def:#1}}
\newcommand{\THM}[1]{Theorem~\ref{thm:#1}}

\newcommand{\SEC}[1]{Section~\ref{sec:#1}}

\DeclareMathOperator{\diam}{diam}
\DeclareMathOperator{\dist}{dist}
\DeclareMathOperator{\USC}{USC}
\DeclareMathOperator{\LSC}{LSC}
\DeclareMathOperator{\osc}{osc}
\DeclareMathOperator{\Lip}{Lip}

\newcommand{\R}{\ensuremath{\mathbb{R}}}
\newcommand{\ep}{\varepsilon}
\newcommand{\mmbox}[1]{\quad \mbox{#1} \quad}
\newcommand{\mrbox}[1]{\quad \mbox{#1} \ }

\begin{document}

\title[An infinity {L}aplace equation with mixed boundary conditions]{An infinity {L}aplace equation with gradient term and mixed boundary conditions}

\author[S. N. Armstrong]{Scott N. Armstrong}
\address{Department of Mathematics, Louisiana State University, Baton Rouge, LA 70803.}
\email{armstrong@math.lsu.edu}
\author[C. K. Smart]{Charles K. Smart}
\address{Department of Mathematics, University of California, Berkeley, CA 94720.}
\email{smart@math.berkeley.edu}
\author[S. J. Somersille]{Stephanie J. Somersille}
\address{Department of Mathematics, University of Texas, Austin TX 78712.}
\email{steph@math.utexas.edu}

\date{\today}


\begin{abstract}
We obtain existence, uniqueness, and stability results for the modified 1-homogeneous infinity Laplace equation \[ -\Delta_\infty u - \beta |Du| = f, \] subject to Dirichlet or mixed Dirichlet-Neumann boundary conditions. Our arguments rely on comparing solutions of the PDE to subsolutions and supersolutions of a certain finite difference approximation.
\end{abstract}

\maketitle

%%%%%
%%%%%
%%%%%

\section{Introduction}

Recently, the first two authors \cite{Armstrong:preprint,Armstrong:inpress} showed that solutions of the 1-homogeneous infinity Laplace equation
\begin{equation}\label{eq:infinity-Laplace}
-\Delta_\infty u := |Du(x)|^{-2} \langle D^2 u(x) Du(x), Du(x) \rangle = f(x)
\end{equation}
can be perturbed by $O(\ep)$ to produce subsolutions and supersolutions of an $\ep$-step finite difference equation. This observation was repeatedly employed in \cite{Armstrong:preprint} to simplify and generalize many aspects of the theory for the PDE \EQ{infinity-Laplace}. 

In this article, we extend this idea to the mixed Dirichlet-Neumann boundary-value problem
\begin{equation}
\label{eq:pde}
\left\{ \begin{aligned}
& - \Delta_\infty u - \beta |Du| = f & \mrbox{in} & \ \Omega, \\
& D_\nu u = 0 & \mrbox{on} & \ \Gamma_N, \\
& u = g & \mrbox{on} & \ \Gamma_D,
\end{aligned} \right.
\end{equation}
where $\Omega \subseteq \R^n$ is a smooth bounded open set, $\Gamma_D \cup \Gamma_N = \partial \Omega$ is a partition of $\partial \Omega$ with $\Gamma_D$ nonempty and closed, $\nu$ is the outer unit normal vector to $\partial \Omega$, and $\beta \in \R$. In the case $\Gamma_N\neq \emptyset$, we also require $\Omega$ to be convex. Our main result, \THM{max-ball}, extends \cite[Proposition 5.3]{Armstrong:preprint} to the PDE in \EQ{pde}.  It states that a subsolution \EQ{pde} becomes a subsolution of a certain finite difference equation after we \emph{max over $\ep$-balls}. We use this result to obtain new existence, uniqueness, and stability results for the boundary-value problem \EQ{pde}.

Peres, Schramm, Sheffield, and Wilson \cite{Peres:2009} showed that the infinity Laplace equation describes the continuum limit of the value functions of a two-player, random-turn game called \emph{$\ep$-step tug-of-war}. A \emph{biased} version of $\ep$-step tug-of-war, in which the randomness favors one of the players, gives rise to the PDE in \EQ{pde} with $\beta \neq 0$, as shown by Peres, G\'abor, and the third author \cite{Peres:preprint} in the case $f \equiv 0$ and $\Gamma_N = \emptyset$. Probabilistic arguments have also been used by Charro, Garc{\'{\i}}a Azorero, and Rossi \cite{Charro:2009,Charro:preprint} to obtain existence of solutions to the infinity Laplace equation with mixed boundary conditions. While we do not use any probabilistic arguments in this paper, the finite difference equation we introduce below does correspond to a certain biased $\ep$-step tug-of-war game. Thus, many of our techniques and arguments have probabilistic analogues. However, the finite difference equation we consider here is designed with analytic consequences in mind (particularly \THM{max-ball}).

In each of the articles mentioned in the paragraph above, solutions of the infinity Laplace equation are obtained by computing the limit as $\ep \to 0$ of the value functions of an $\ep$-step tug-of-war game. In contrast, in this paper we first establish an analogue of \emph{comparisons with cones}, which allows us to apply the Perron process to establish the existence of maximal and minimal solutions of our boundary-value problem. This approach to existence has also been employed recently by Lu and Wang \cite{Lu:preprint} in the case $\beta =0$, $\Gamma_N=\emptyset$. Our uniqueness and stability results rely on the perturbation theorem (\THM{max-ball}), which allows us to reduce our consideration to the finite difference equation.

In the next section, we review the notion of viscosity solution for the mixed boundary-value problem \EQ{pde} and then state our assumptions and main results. In \SEC{perturbation} we prove the perturbation theorem. In \SEC{etc} we prove existence and apply the perturbation theorem to obtain uniqueness and stability results for boundary-value problem.

\section{Statement of main results} \label{sec:prelim}

We take $\Omega \subseteq \R^n$ to be a bounded open connected domain with a $C^1$ boundary $\partial \Omega$. We partition the boundary $\partial \Omega$ into subsets $\Gamma_N = \partial \Omega \cap A$ and $\Gamma_D = \partial \Omega \setminus A \neq \emptyset$, for some open subset $A \subseteq \R^n$. In the case $\Gamma_N \neq \emptyset$, we require the domain $\Omega$ to be convex. We denote the outward pointing unit normal vector to $\partial \Omega$ by $\nu$.

Let us recall the notion of viscosity solution of \EQ{pde}. Define for $\varphi \in C^2(\R^n)$ the operators
\begin{equation*}
\Delta_\infty^+ u(x) := \left\{ \begin{array}{ll}
|D\varphi(x)|^{-2} \langle D^2 \varphi(x) D\varphi(x), D\varphi(x) \rangle & \mrbox{if} D\varphi(x) \neq 0, \\
\max \left\{ \langle D^2 \varphi (x) v, v \rangle : |v| =  1 \right\} & \mrbox{if} D\varphi(x) = 0,
\end{array} \right.
\end{equation*}
and $\Delta_\infty^- \varphi(x) := - \Delta_\infty^+(-\varphi)(x)$. We denote the set of real-valued upper semicontinuous functions on a set $V \subseteq \R^n$ by $\USC(V)$, and the set of real-valued lower semicontinuous functions by $\LSC(V)$. The set of Lipschitz functions on $V$ is denoted by $\Lip(V)$, and the Lipschitz constant of $u\in \Lip(V)$ is denoted by $\Lip(u,V)$. The oscillation of a function $u : V \to \R$ is denoted $\osc u := \sup_V u - \inf_V u$.

\begin{definition}\label{def:visc}
Given $h\in \USC(\Omega \cup \Gamma_N)$, we say that $u \in \USC(\bar \Omega)$ is a \emph{viscosity subsolution} of the system
\begin{equation}
\label{eq:PDE-Neumann}
\left\{ \begin{aligned}
& - \Delta_\infty u - \beta |Du| = h & \mrbox{in} & \ \Omega, \\
& D_\nu u = 0 & \mrbox{on} & \  \Gamma_N,
\end{aligned} \right.
\end{equation}
if, for every $\varphi \in C^2(\bar \Omega)$ and $x_0 \in \Omega \cup \Gamma_N$ such that the map $x \mapsto (u - \varphi)(x)$ has a local maximum at $x_0$, we have
\begin{equation*}
- \Delta_\infty^+ \varphi(x_0) - \beta |D\varphi(x_0)| \leq h(x_0)
\end{equation*}
or
\begin{equation*}
x_0 \in \Gamma_N \quad \mbox{and} \quad D\varphi (x_0) \cdot \nu(x_0) \leq 0.
\end{equation*}
Similarly, if $h\in \LSC(\Omega \cup \Gamma_N)$, then $v \in \LSC(\bar \Omega)$ is a \emph{viscosity supersolution} of \EQ{PDE-Neumann} if $-v$ is a subsolution of \EQ{PDE-Neumann} for $-h$. If $h \in C(\Omega \cup \Gamma_N)$, then $u \in C(\bar \Omega)$ is a \emph{viscosity solution} of \EQ{PDE-Neumann} if it is both a viscosity subsolution and a viscosity supersolution of \EQ{PDE-Neumann} for $h$.
\end{definition}

We remark that by altering \DEF{visc} by requiring all local extrema to be strict, we obtain an equivalent definition of viscosity subsolution and viscosity supersolution. We emphasize that all differential inequalities in this paper involving functions not known to be smooth are to be understood in the viscosity sense. In particular, if we say that $u$ is a solution of the system of differential inequalities
\begin{equation*}
\left\{ \begin{aligned}
& - \Delta_\infty u - \beta |Du| \leq h & \mrbox{in} & \ \Omega, \\
& D_\nu u \leq 0 & \mrbox{on} & \  \Gamma_N,
\end{aligned} \right.
\end{equation*}
this is taken to mean that $u$ is a viscosity subsolution of \EQ{PDE-Neumann}.

For each $\ep > 0$, we define $\Omega_\ep$ to be the set of points in $\bar \Omega$ which are farther than $\ep$ from the Dirichlet boundary $\Gamma_D$. That is,
\begin{equation*}
\Omega_\ep := \{ x \in \bar \Omega : \dist(x, \Gamma_D) > \ep \}.
\end{equation*}
If $h, -\tilde h \in \USC(\Omega \cup \Gamma_N)$ and $\ep > 0$, we denote
\begin{equation*}
h^\ep(x) := \max_{\bar B(x,\ep)} h \mmbox{and} \tilde h_\ep(x) := \min_{\bar B(x,\ep)} \tilde h, \quad \mbox{for} \  x\in \Omega_\ep.
\end{equation*}
For $u\in C(\Omega \cup \Gamma_N)$, we define the quantities
\begin{equation*}
S^+_\ep u(x) := \frac{1}{\ep} \left( u^\ep(x) - u(x) \right) \mmbox{and} S^-_\ep u(x) := \frac{1}{\ep} \left( u(x) - u_\ep(x) \right).
\end{equation*}
Our finite difference approximation to \EQ{PDE-Neumann} is the scheme
\begin{equation}\label{eq:fd}
a^-_\ep(\beta) S^-_\ep u(x) - a^+_\ep(\beta) S^+_\ep u(x) = \ep h(x),
\end{equation}
where the coefficients $a^\pm_\ep(\beta)$ are given by
\begin{equation*}
a^+_\ep(\beta) := \frac{\beta}{\exp(\ep\beta) - 1} \quad \mbox{and} \quad a^-_\ep(\beta) := \frac{\beta}{1-\exp(-\ep\beta)}, \quad \beta \in \R \setminus \{ 0 \},
\end{equation*}
and $a^+_\ep(0) := a^-_\ep(0) := 1/\ep$.

\medskip

Our main result is the following perturbation theorem.

\begin{thm}
\label{thm:max-ball}
Suppose that $h \in \USC(\Omega \cup \Gamma_N)$ and $u \in C(\Omega \cup \Gamma_N)$ is viscosity subsolution of the system \EQ{PDE-Neumann}. Then for each $\ep > 0$,
\begin{equation}
\label{eq:fd-subsol}
a^-_\ep(\beta) S^-_\ep u^\ep(x) - a^+_\ep(\beta) S^+_\ep u^\ep(x) \leq h^{2 \ep}(x) \mmbox{for every} x \in \Omega_{2\ep}.
\end{equation}
\end{thm}

\THM{max-ball} asserts that small perturbations of viscosity subsolutions of \EQ{pde} are subsolutions of the finite difference equation \EQ{fd}. We will see below that \EQ{fd} has a simple comparison lemma (see \LEM{fd-comparison}), which we combine with \THM{max-ball} to deduce the following comparison result for \EQ{PDE-Neumann}.

\begin{thm}
\label{thm:comparison}
Assume that $h, - \tilde h \in \USC(\Omega \cup \Gamma_N) \cap L^\infty(\Omega\cup \Gamma_N)$ satisfy $h\leq \tilde h$, as well as
\begin{equation*}
h < \tilde h, \quad h\equiv 0, \quad \mbox{or} \quad \tilde h > 0.
\end{equation*}
Suppose that $u \in \USC(\bar \Omega)$ is a viscosity subsolution of the system
\begin{equation}\label{eq:compare-1}
\left\{ \begin{aligned}
& - \Delta_\infty u - \beta |Du| \leq h & \mrbox{in} & \ \Omega, \\
& D_\nu u \leq 0 & \mrbox{on} & \ \Gamma_N,
\end{aligned} \right.
\end{equation}
and $v \in \LSC(\bar \Omega)$ is a viscosity supersolution of the system
\begin{equation}\label{eq:compare-2}
\left\{ \begin{aligned}
& - \Delta_\infty v - \beta |Dv| \geq \tilde h & \mrbox{in} &  \  \Omega, \\
& D_\nu v \geq 0 & \mrbox{on} & \ \Gamma_N.
\end{aligned} \right.
\end{equation}
Then
\begin{equation}
\label{eq:comparison-conclusion}
\max_{\bar \Omega} (u - v) = \max_{\Gamma_D} (u - v).
\end{equation}
\end{thm}

We will establish the following existence result for \EQ{pde} using the Perron method, which extends \cite[Theorem 2.14]{Armstrong:preprint}.

\begin{thm} \label{thm:existence}
Assume that $f \in C(\Omega \cup \Gamma_N) \cap L^\infty( \Omega \cup \Gamma_N)$ and $g \in C(\Gamma_D)$. Then there exist solutions $\underline u, \overline u \in C(\bar \Omega)$ of the boundary-value problem \EQ{pde}, with the property that $u \leq \overline u$ ($u \geq \underline u$) if $u$ is a subsolution (supersolution) of \EQ{pde} with $u \leq g$ ($u \geq g$) on $\Gamma_D$.
\end{thm}

Using \THM{max-ball} we obtain the following stability result.

\begin{thm}
\label{thm:stability}
Suppose that $\beta, \beta_j \in \R$, $f, f_j \in \USC(\Omega \cup \Gamma_N)$, and $u, u_j \in C( \Omega \cup \Gamma_N)$ such that for each $j\geq 1$, the function $u_j$ satisfies the system
\begin{equation*}
\left\{ \begin{aligned}
& -\Delta_\infty u_j - \beta_j |Du_j| \leq f_j & \mbox{in} & \ \Omega, \\
& D_\nu u_j \leq 0 & \mbox{on} & \ \Gamma_N.
\end{aligned} \right.
\end{equation*}
Suppose in addition that $\beta_j \to \beta$, and $f_j \to f$ and $u_j \to u$ locally uniformly in $\Omega\cup\Gamma_N$. Then $u$ is a subsolution of the system
\begin{equation*}
\left\{ \begin{aligned}
& -\Delta_\infty u - \beta |Du| \leq f & \mbox{in} & \ \Omega, \\
& D_\nu u \leq 0 & \mbox{on} & \ \Gamma_N.
\end{aligned} \right.
\end{equation*}
\end{thm}

Combining sup-norm and interior Lipschitz estimates (see \LEM{local-lipschitz}, below) with \THM{stability}, we  immediately obtain the following result.

\begin{cor}
Let $\beta$, $\beta_j$, $f$, $f_j$ be as in the hypotheses of \THM{stability}, and suppose in addition that $|f_j| \leq K$ for all $j$. Assume that $g, g_j\in C(\Gamma_D)$ such that $g_j \to g$ uniformly, and $u_j$ is a solution of the problem
\begin{equation*}
\left\{ \begin{aligned}
& -\Delta_\infty u_j - \beta_j |Du_j| = f_j & \mbox{in} & \ \Omega, \\
& D_\nu u_j = 0 & \mbox{on} & \ \Gamma_N,\\
& u_j = g_j & \mbox{on} & \ \Gamma_D.
\end{aligned} \right.
\end{equation*}
Then there is a subsequence $\{u_{j_k}\}$ and a solution $u\in C(\bar \Omega)$ of the problem
\begin{equation*}
\left\{ \begin{aligned}
& -\Delta_\infty u - \beta |Du| = f & \mbox{in} & \ \Omega, \\
& D_\nu u = 0 & \mbox{on} & \ \Gamma_N,\\
& u = g & \mbox{on} & \ \Gamma_D,
\end{aligned} \right.
\end{equation*}
such that $u_{j_k} \to u$ uniformly on $\bar \Omega$.
\end{cor}

%%%%%
%%%%%
%%%%%

\section{The perturbation theorem} \label{sec:perturbation}

We begin our proof of \THM{max-ball} with an analogue of \emph{comparisons with cones from above} for infinity subharmonic functions (see \cite{Crandall:2001}).

\begin{lem}
\label{lem:cca}
Suppose $k \in \R$ and $u \in \USC(\bar \Omega)$ is a subsolution of
\begin{equation}
\left\{ \begin{aligned}
\label{eq:subsol-const}
& - \Delta_\infty u - \beta |Du| \leq k & \mrbox{in} & \ \Omega, \\
& D_\nu u \leq 0 & \mrbox{on} & \ \Gamma_N.
\end{aligned} \right.
\end{equation}
Suppose further that $r>0$ and $\gamma \in C^2([0, r])$ satisfies
\begin{equation*}
- \gamma'' - \beta |\gamma'| = k \mmbox{and} \gamma' > 0 \mrbox{in} (0,r),
\end{equation*}
and define
\begin{equation*}
\varphi(x) := \gamma(|x - x_0|)
\end{equation*}
for some $x_0 \in \bar \Omega$. Set $\Lambda :=(\partial B(x_0,r) \cap \bar\Omega) \cup (B(x_0,r) \cap \Gamma_D)$. Then
\begin{equation}\label{eq:gamma-p-neq}
\max_{\bar B(x_0,r) \cap \bar \Omega} (u - \varphi) = \max_{\Lambda \cup \{ x_0 \}} (u - \varphi),
\end{equation}
and if in addition $\gamma'(0) = 0$, then
\begin{equation} \label{eq:gamma-p-eq}
\max_{\bar B(x_0,r) \cap \bar \Omega} (u - \varphi) = \max_{\Lambda} (u - \varphi).
\end{equation}
\end{lem}

\begin{proof}
By the upper semicontinuity of $u$, we may assume that $x_0 \in \Omega$. By the convexity of $\Omega$ in the case that $\Gamma_N \neq \emptyset$, we have that $(y - x_0) \cdot \nu > 0$ for all $y \in \Gamma_N$.

Consider first the case that $\gamma'(0) >0$.  Arguing indirectly, suppose that the conclusion \EQ{gamma-p-neq} fails so that we can find $x_1 \in \bar B(x_0,r) \cap \bar \Omega$ such that
\begin{equation}\label{eq:point-max}
(u-\varphi)(x_1) > \max_{\Lambda \cup \{ x_0 \}} (u - \varphi).
\end{equation}
For $\delta, \eta \geq 0$ to be chosen below, let $\tilde \gamma \in C^2([0,r])$ be the solution of the initial-value problem
\begin{equation*}
- \tilde \gamma'' - \beta  |\tilde \gamma'| = k + \delta, \quad \tilde \gamma(0) = \gamma(0), \quad \tilde \gamma'(0) = \gamma'(0) + \eta.
\end{equation*}
For the choice $\eta > 0$ and $\delta =0$, it is easy to check that $\tilde \gamma' > \gamma' \geq 0$ on $[0,r]$. Thus by elementary stability properties of our ODE, we can select small enough $\eta > 0$ and $\delta > 0$ such that $\tilde\gamma' > 0$ on $[0,r]$ and the function $\tilde \varphi(x) : = \tilde\gamma(|x-x_0|)$ satisfies \EQ{point-max} with $\varphi$ replaced by $\tilde \varphi$.

Select a point $x_2 \in \bar B(x_0, r) \cap \bar \Omega$ at which the map $x \mapsto (u - \tilde \varphi)(x)$ has attains its maximum in $\bar B(x_0, r) \cap \bar \Omega$, and note that $x_2 \not\in \Lambda \cup \{ x_0 \}$. Observe that $\tilde \varphi$ is $C^2$ in a neighborhood of $x_2$ and
\begin{equation*}
- \Delta_\infty^- \tilde \varphi(x_2) - \beta |D \tilde \varphi(x_2)| = k + \delta > k.
\end{equation*}
Since $u$ is a viscosity subsolution of \EQ{subsol-const}, it must be the case that $x_2 \not \in \Omega$. Thus  we must have $x_2 \in \Gamma_N$. Using again that $u$ is a viscosity subsolution of \EQ{subsol-const}, we have $\nu \cdot D\varphi(x_2)  \leq 0$. This contradicts the fact that $\tilde \gamma'(|x_2-x_0|)  > 0$ and $(x_2 - x_0) \cdot \nu > 0$. We have verified \EQ{gamma-p-neq} in the case $\gamma'(0) > 0$.

In the case $\gamma'(0) = 0$ we also argue indirectly, and suppose that \EQ{gamma-p-eq} fails, so that there exists $x_1 \in \bar B(x_0,r) \cap \bar \Omega$ such that 
\begin{equation}\label{eq:point-max-0}
(u-\varphi)(x_1) > \max_{\Lambda} (u - \varphi).
\end{equation}
It is easy to see that $\gamma'(0) = 0$ and $\gamma'>0$ on $(0,r)$ implies that $k< 0$. For $0 \leq \delta \leq -k$ to be selected below, let $\tilde \gamma$ be defined by
\begin{equation*}
\tilde\gamma(t) := \begin{cases} -\frac{k+\delta}{\beta} t - \frac{k+\delta}{\beta^2} \left( e^{-\beta t} -1 \right) +\gamma(0), & \mbox{if} \ \beta \neq 0, \\
\frac{k+\delta}{2}t^2 + \gamma(0), & \mbox{if} \ \beta =0.
\end{cases}
\end{equation*}
Then $\tilde \gamma \in C^2([0,\infty))$ satisfies $\tilde \gamma' > 0$ in $(0,\infty)$ and
\begin{equation*}
- \tilde \gamma'' - \beta  |\tilde \gamma'| = k + \delta, \quad \tilde \gamma(0) = \gamma(0), \quad \tilde \gamma'(0) = 0.
\end{equation*}
For sufficiently small $0< \delta < -k$, the function $\tilde \varphi(x):= \tilde \gamma(|x-x_0|)$ satisfies \EQ{point-max-0} with $\varphi$ replaced by $\tilde \varphi$. Notice also that $\tilde\varphi\in C^2(\R^n)$, and
\begin{equation*}
-\Delta^-_\infty \tilde\varphi(x) -\beta |D\tilde\varphi(x)| = k+\delta > k, \quad x\in  \R^n.
\end{equation*}
Select a point $x_2 \in \bar B(x_0, r) \cap \bar \Omega$ at which the map $x \mapsto (u - \tilde \varphi)(x)$ attains its maximum in $\bar B(x_0,r) \cap \bar \Omega$, and note that $x_2\not \in \Lambda$. Recalling that $x_0 \not\in \partial \Omega$, we may proceed as above to derive a contradiction.
\end{proof}

The following lemma is a version of \THM{max-ball} for our cone functions.

\begin{lem}
\label{lem:max-ball-1d}
Assume that $\gamma \in C^2([0, r])$ satisfies
\begin{equation*}
- \gamma'' - \beta |\gamma'| = k \mmbox{and} \gamma' > 0 \mrbox{in} (0,r).
\end{equation*}
Suppose that $\ep > 0$ and $r \geq 2 \ep$ or $\gamma'(r) = 0$. Then
\begin{equation*}
a^+_\ep(\beta) (\gamma(r_1) - \gamma(0)) - a^-_\ep(\beta) (\gamma(r_2) - \gamma(r_1) \leq \ep k,
\end{equation*}
where $r_1 := \min\{\ep, r\}$ and $r_2 := \min\{2\ep, r\}$.
\end{lem}

\begin{proof}
We consider only the case $\beta \neq 0$. The case $\beta = 0$ is similar (and easier) and was proved in \cite{Armstrong:preprint}.

Since $\gamma' > 0$ in $(0,r)$, we see that $\gamma$ satisfies the linear equation
\begin{equation*}
- \gamma'' - \beta \gamma' = k \mrbox{in} (0,r).
\end{equation*}
In particular, there are constants $c_1, c_2 \in \R$ such that
\begin{equation} \label{eq:gammac2}
\gamma(t) = c_1 + c_2 e^{- \beta t} - \frac{k}{\beta} t, \quad 0 \leq t \leq r.
\end{equation}
Consider the case $r \geq 2 \ep$. Then we have
\begin{equation*}
\gamma(r_1) - \gamma(0) = c_2 e^{- \beta \ep} (1- e^{\beta \ep} ) - \frac{k}{\beta} \ep,
\end{equation*}
and
\begin{equation*}
\gamma(r_2) - \gamma(r_1) = c_2 e^{- \beta \ep} (e^{- \beta \ep} -1 ) - \frac{k}{\beta} \ep.
\end{equation*}
The conclusion now follows from some algebra, making use of the identity
\begin{equation}
\label{eq:ode-expcalc}
(e^{\beta \ep} - 1)^{-1} - (1 - e^{- \beta \ep})^{-1} = -1.
\end{equation}

Now consider the case $r < 2 \ep$ and $\gamma'(r) = 0$. Note that $\gamma'(r) = 0$ implies
\begin{equation*}
\gamma(t) = c_1 - \frac{k}{\beta ^2} e^{- \beta(t - r)} - \frac{k}{\beta}t.
\end{equation*}
Note also that $\gamma' > 0$ in $(0,r)$ and $\beta \neq 0$ imply $k > 0$. In particular, if we extend $\gamma$ to all of $\R$ by the formula above, we see that it is concave and achieves its maximum at $t=r$. Therefore, using \EQ{gammac2} we obtain
\begin{equation*}
\gamma(r_1) - \gamma(0) \leq \gamma(r_1) - \gamma(r_1 - \ep) = c_2 e^{- \beta r_1} (1-e^{\beta \ep} ) - \frac{k}{\beta} \ep
\end{equation*}
and
\begin{equation*}
\gamma(r_2) - \gamma(r_1) \geq \gamma(r_1+\ep) - \gamma(r_1) = c_2 e^{- \beta r_1} (e^{- \beta \ep}-1) - \frac{k}{\beta} \ep.
\end{equation*}
The conclusion follows from a little algebra, as above.
\end{proof}

\THM{max-ball} follows at once from the following lemma.

\begin{lem}
\label{lem:max-ball}
Suppose that $k \in \R$ and $u \in C(\bar \Omega)$ is a subsolution of \EQ{subsol-const}. Then
\begin{equation*}
a^-_\ep(\beta) S^-_\ep u^\ep(x) - a^+_\ep(\beta) S^+_\ep u^\ep(x)  \leq k \mrbox{for every} x\in \Omega_{2\ep}.
\end{equation*}
\end{lem}

\begin{proof}
Suppose $u \in C(\bar \Omega)$ satisfies \EQ{subsol-const} and $x \in \Omega_{2 \ep}$. Let $y \in \bar B(x,\ep) \cap \bar \Omega$ and $z \in \bar B(x, 2 \ep) \cap \bar \Omega$ satisfy
\begin{equation*}
u(y) = u^\ep(x) \mmbox{and} u(z) = u^{2\ep}(x),
\end{equation*}
and observe that
\begin{equation}\label{eq:1d-to-nd}
\ep S^-_\ep u^\ep(x) \leq u(y) - u(x) \mmbox{and} \ep S^+_\ep u^\ep(x) = u(z) - u(y).
\end{equation}

Select $\delta > 0$, and let $\gamma \in C^2([0,r])$ be the solution of the initial-value problem
\begin{equation*}
- \gamma'' - \beta |\gamma'| = k, \quad \gamma(0) = u(x), \quad \gamma'(0) = \eta,
\end{equation*}
where $\eta$ is minimized subject to the constraints $\eta \geq 0$ and $\max_{[0,2\ep]} \gamma \geq u(z) + \delta$. Set $r: = \inf\{ s > 0 : \gamma(s) \geq u(z)+ \delta \}$, and notice that $\gamma(r) = u(z) + \delta > u(x) = \gamma(0)$, and thus $0 < r \leq 2\ep$.

We claim that $\eta > 0$. Suppose on the contrary that $\eta =0$. Since $\gamma$ is not constant, it is easy to see that its derivative $\gamma'(t)$ can only vanish at a single point $t\in \R$. Since $\gamma'(0) =0$ and $\gamma (r) > \gamma(0)$, we deduce that $\gamma'> 0$ on $(0,\infty)$. Define $\varphi(z) := \gamma(|z-x|) - \delta$, and notice that $\varphi \geq u(z) \geq u$ on $\partial B(x,2\ep) \cap \bar \Omega$. From \LEM{cca} we deduce that $\varphi \geq u$ in $B(x,2\ep) \cap \bar \Omega$, which is not possible since $\varphi(x) = u(x) -\delta < u(x)$. Thus $\eta > 0$. It follows that $\max_{[0,2\ep]} \gamma = u(z)+\delta$.

Next we show that
\begin{equation}\label{eq:weirdo-claim}
\gamma'> 0 \quad \mbox{in} \ (0,r) \quad \mbox{and} \quad r=2\ep \ \ \mbox{or} \ \ \gamma'(r) = 0.
\end{equation}
As mentioned above, since $\gamma$ is not constant its derivative $\gamma'(t)$ can only vanish at a single point $t\in \R$. Moreover, at this point $\gamma'$ must change sign. Owing to the initial condition $\gamma'(0) = \eta > 0$, we see that if $\gamma'(t) = 0$ for some $t>0$, then $\gamma$ achieves its maximum at $t$. Since $\max_{[0,2\ep]} \gamma = u(z) + \delta = \gamma(r)$, the claim \EQ{weirdo-claim} follows.

Define $\tilde \varphi(w):= \gamma(|w-x|)$. Since $\tilde\varphi \geq u$ on $\{ x \} \cup \partial B(x, r) \cap \bar \Omega$, \LEM{cca} implies that $\tilde\varphi \geq u$ on $\bar B(x,r) \cap \bar \Omega$. In particular, $u(y) \leq \gamma(\min\{ r, \ep\})$. Since $u(z) + \delta = \gamma(r) = \gamma(\min\{ r, 2 \ep \})$, the conclusion follows from \LEM{max-ball-1d} and \EQ{1d-to-nd} once we send $\delta \to 0$.
\end{proof}

\begin{proof}[{\bf Proof of \THM{max-ball}}]
Our result is immediately obtained from \LEM{max-ball}.
\end{proof}

%%%%%
%%%%%
%%%%%

\section{Existence, uniqueness, and stability} \label{sec:etc}

In this section we prove the rest of our main results, Theorems \ref{thm:comparison}, \ref{thm:existence}, and \ref{thm:stability}. First, we need a sup-norm and Lipschitz estimate for subsolutions of \EQ{PDE-Neumann}, which we obtain with the help of \LEM{cca}.

\begin{lem}
\label{lem:local-lipschitz}
Suppose that $k\geq 0$ and $u \in \USC(\bar \Omega)$ is a subsolution of
\begin{equation*}
\left\{ \begin{aligned}
& -\Delta_\infty u - \beta |Du| \leq k & \mbox{in} & \ \Omega, \\
& D_\nu u \leq 0 & \mbox{on} & \ \Gamma_N.
\end{aligned} \right.
\end{equation*}
Then 
\begin{equation} \label{eq:u-sup}
\max_{\bar\Omega} u  \leq \max_{\Gamma_D} u + C_1,
\end{equation}
where the constant $C_1>0$ depends only on $k$, $\beta$, and $\diam(\Omega)$. If in addition $u$ is bounded below, then $u\in \Lip(\Omega_\delta)$ for each $\delta > 0$, and we have the estimate
\begin{equation*}
| u(x) - u(y) | \leq C_2 |x-y| \quad \mbox{for every} \ \ x,y\in \Omega_\delta,
\end{equation*}
where the constant $C_2$ depends only on $\osc u$, $\delta$, $k$, $\beta$, and $\diam(\Omega)$.
\end{lem}

\begin{proof}
Denote $d:= \diam(\Omega)$ and define a function $\gamma=\gamma(t)$ by
\begin{equation*}
\gamma(t) : = \begin{cases}
-\frac{1}{\beta} k t - \frac{1}{\beta^2}k e^{\beta d} \left(e^{-\beta t} -1\right) & \mbox{if} \ \beta \neq 0, \\
-\frac{1}{2}k t^2 + dk t & \mbox{if} \ \beta = 0, \\
\end{cases}
\end{equation*}
Observe that $\gamma$ satisfies $\gamma(0) = 0$,
\begin{equation*}
- \gamma'' - \beta |\gamma'| = k \quad \mbox{in} \ (0,d), \quad \mbox{and} \quad \gamma' > 0 \quad \mbox{in} \ (0,d).
\end{equation*}
Select $x_0 \in \Gamma_D$, and set
\begin{equation} \label{eq:smooth-supersolution}
\varphi(x) := \max_{\Gamma_D} u + \gamma(|x - x_0|).
\end{equation}
By \LEM{cca}, $u \leq \varphi$ on $\bar \Omega$. In particular,
\begin{equation*}
u \leq \max_{\bar \Omega} \varphi \leq  \max_{\Gamma_D} u  + \gamma(d).
\end{equation*}
Thus we have the estimate \EQ{u-sup} for $C_1:= \gamma(d)$.

If $u$ is bounded below and $\delta > 0$, then for any $x\in \Omega_\delta$ and $y\in B(x,\delta) \cap \bar\Omega$,
\begin{equation*}
u(y) \leq \varphi^x(y) := u(x) + \max\{ 1, \osc u / \gamma(\delta) \} \gamma(|y-x|).
\end{equation*}
This implies that $\Lip(u,\Omega_\delta) \leq \max\{ 1, \osc u / \gamma(\delta) \} \gamma'(0)=:C_2$.
\end{proof}

We now establish a comparison result for solutions of our finite difference equation. As in \cite{Armstrong:preprint,Armstrong:inpress}, our simple argument uses an idea of Le Gruyer \cite{LeGruyer:2007}.

\begin{lem}
\label{lem:fd-comparison}
Suppose $u, v \in C(\bar \Omega)$ and $f, \tilde f :\Omega_\ep \to \R$ satisfy
\begin{equation} \label{eq:FDeq}
a^-_\ep(\beta) S^-_\ep u - a^+_\ep(\beta)S^+_\ep u \leq f \leq \tilde f \leq a^-_\ep(\beta) S^-_\ep v - a^+_\ep(\beta)S^+_\ep v \mrbox{in} \Omega_\ep.
\end{equation}
Suppose in addition that (i) $f < \tilde f$, (ii) $f \leq 0$, or (iii) $\tilde f \geq 0$. Then
\begin{equation*}
\max_{\bar \Omega} (u - v) = \max_{\bar \Omega \setminus \Omega_\ep} (u - v).
\end{equation*}
\end{lem}

\begin{proof}
Arguing indirectly, let us suppose that the hypothesis holds but the conclusion fails. Define
\begin{equation*}
E := \{ x \in \bar \Omega : (u - v)(x) = \max_{\bar \Omega} (u - v) \},
\end{equation*}
and observe that $E$ is nonempty, closed, and contained in $\Omega_\ep$. Select $x_0 \in E$. Since the map $x \mapsto (u - v)(x)$ attains its maximum at $x_0$, we have
\begin{equation*}
S^-_\ep u(x_0) \geq S^-_\ep v(x_0) \mmbox{and} S^+_\ep u(x_0) \leq S^+_\ep v(x_0).
\end{equation*}
Since $f(x_0) \leq \tilde f(x_0)$, and $a^\pm_\ep(\beta) > 0$, the inequalities in \EQ{FDeq} must be equalities at $x_0$, and thus $f(x_0) = \tilde f(x_0)$,
\begin{equation}
\label{eq:slopes-equal}
S^-_\ep u(x_0) = S^-_\ep v(x_0) \mmbox{and} S^+_\ep u(x_0) = S^+_\ep v(x_0).
\end{equation}
In particular, in the case (i) holds we obtain a contradiction. 

We have left to consider the cases (ii) and (iii). By symmetry, we may only consider (ii). Define
\begin{equation*}
F := \{ x \in E : u(x) = \max_E u \},
\end{equation*}
and notice that $F$ is nonempty, closed, and $F\subseteq E \subseteq \Omega_\ep$.

Suppose there is an $x_0 \in F$ and an $x_1 \in \bar B(x_0, \ep)$ such that
\begin{equation*}
u(x_1) - u(x_0) = \ep S^+_\ep u(x_0) > 0.
\end{equation*}
Since $x_1 \notin E$, we must have $u(x_1) - v(x_1) < u(x_0) - v(x_0)$. In particular
\begin{equation*}
\ep S^+_\ep v(x_0) \geq v(x_1) - v(x_0) > u(x_1) - u(x_0) = \ep S^+ u(x_0),
\end{equation*}
contradicting \EQ{slopes-equal}. Thus $S^+_\ep u \equiv 0$ on $F$.

Since $f \leq 0$, by \EQ{FDeq} we must have $S^-_\ep u \equiv 0$ on $F$. Thus $u$ is constant on $\bar B(x_0, \ep)$ for all $x_0 \in F$. By \EQ{slopes-equal}, the same is true for $v$. As $F\subseteq \Omega_\ep$ is closed, we can choose a point $x_0\in \partial F \subseteq F$ to obtain a contradiction.
\end{proof}

\begin{proof}[\bf Proof of \THM{comparison}]
We first reduce to the case that $u, v\in C(\Omega \cup \Gamma_N)$. It is clear from \DEF{visc} that the minimum of two supersolutions is also a supersolution. Observe that the function $\varphi$ defined in \EQ{smooth-supersolution} is a smooth supersolution of \EQ{compare-2}, where we take $k:= \| \tilde h \|_{L^\infty(\Omega \cup \Gamma_N)}$. Thus $\tilde v:= \min\{ v , \varphi \}$ is a supersolution of \EQ{compare-2} which is bounded above. By \LEM{local-lipschitz}, we see that $\tilde v\in C(\Omega \cup \Gamma_N)$. By replacing $v$ with $\tilde v$, we may assume that $v \in C(\Omega \cup \Gamma_N)$. Similarly, we may assume that $u \in C(\Omega \cup \Gamma_N)$.

According to \THM{max-ball}, we have
\begin{equation*}
a^-_\ep(\beta) S^-_\ep u^\ep - a^+_\ep(\beta) S^+_\ep u^\ep \leq \ep h^{2 \ep} \quad \mbox{and} \quad a^-_\ep(\beta) S^-_\ep v_\ep - a^+_\ep(\beta) S^+_\ep v_\ep \geq \ep \tilde h_{2 \ep} \quad \mbox{in}\  \Omega_{2\ep}.
\end{equation*}
In the case $h \equiv 0 \leq \tilde h$, we have $h^{2\ep} \equiv 0 \leq \tilde h_{2\ep}$ and so may apply \LEM{fd-comparison} to deduce that
\begin{equation}
\label{eq:ep-comparison}
\max_{\bar \Omega_{\ep}} (u^\ep - v_\ep) = \max_{\bar \Omega_{\ep} \setminus \Omega_{2\ep}} (u^\ep - v_\ep).
\end{equation}
In the case $h < \tilde h$, for each $r > 0$ we may choose $0 < \ep < r/2$ such that $h^{2\ep} < \tilde h_{2 \ep}$ in $\Omega_r$, and then apply \LEM{fd-comparison} to obtain
\begin{equation}
\label{eq:ep-r-comparison}
\max_{\bar \Omega_{r - \ep}} (u^\ep - v_\ep) = \max_{\bar \Omega_{r - \ep} \setminus \Omega_{r}} (u^\ep - v_\ep).
\end{equation}
Sending $\ep \to 0$ in \EQ{ep-comparison}, and $\ep \to 0$ followed by $r \to 0$ in \EQ{ep-r-comparison}, and using the upper semicontinuity of $u$ and $-v$ up to the boundary, we obtain \EQ{comparison-conclusion} in the cases $h< \tilde h$ and $h\equiv 0$.

We have left to consider the case $\tilde h > 0$. Define $w:= (1+\delta) v$ for $\delta >0$, and observe that $w$ satisfies 
\begin{equation*}
\left\{ \begin{aligned}
& -\Delta_\infty w - \beta |Dw| \geq (1+\delta) \tilde h & \mbox{in} & \ \Omega, \\
& D_\nu w = 0 & \mbox{on} & \ \Gamma_N.
\end{aligned} \right.
\end{equation*}
Since $(1+\delta) \tilde h > h$, our results above imply
\begin{equation*}
\max_{\bar \Omega} (u-(1+\delta) v) = \max_{\Gamma_D} (u-(1+\delta)v).
\end{equation*}
Sending $\delta \to 0$, we obtain \EQ{comparison-conclusion}.
\end{proof}

\begin{proof}[\bf Proof of \THM{existence}]
Define $u : \bar \Omega \rightarrow \R$ by
\begin{equation*}
u(x) := \sup \{ w(x) : w\in \USC(\bar\Omega) \ \mbox{is a subsolution of \EQ{pde} and} \ w \leq g \ \mbox{on} \ \Gamma_D  \}.
\end{equation*}
To see that $u$ is well-defined, denote $d:= \diam(\Omega)$, $k:= \| f \|_{L^\infty(\Omega \cup\Gamma_N)}$, and set
\begin{equation*}
\gamma(t) : = \begin{cases}
-\frac{1}{\beta} k t + \frac{1}{\beta^2}k e^{-\beta d} \left(e^{\beta t} -1\right) & \mbox{if} \ \beta \neq 0, \\
\frac{1}{2}k t^2 - dk t & \mbox{if} \ \beta = 0. \\
\end{cases}
\end{equation*}
Observe that $\gamma(0)=0$, $\gamma'< 0$ on $(0,d)$, and 
\begin{equation*}
-\gamma'' - \beta |\gamma'| = -k\quad \mbox{in} \ (0,d).
\end{equation*}
Select $x_0 \in \Gamma_D$ and define $\varphi(x) := \min_{\Gamma_D} g + \gamma(|x-x_0|)$. It is clear that $\varphi$ is a smooth subsolution of \EQ{pde} and $\varphi \leq g$ on $\Gamma_D$. Thus $u \geq \varphi$, and in particular $u$ is bounded below. According to \LEM{local-lipschitz}, $u$ is also bounded above. By construction, $u \in \USC(\bar\Omega)$. 

We now proceed to show that $u=g$ on $\Gamma_D$. By construction, $u\leq g$ on $\Gamma_D$. To get the other inequality, select $y \in \Gamma_D$, and $\ep > 0$. For $r> 0$ small enough, $g\geq g(y) - \ep$ on $\Gamma_D \cap \bar B(y,r)$. With $\gamma$ as above, define
\begin{equation*}
\varphi^y(x):= g(y) - \ep + \max\left\{ 1, \frac{\min_{\Gamma_D} g - g(y)}{\gamma(r)} \right\} \gamma(|x-y|).
\end{equation*}
Observe that $\varphi^y$ is a smooth subsolution of \EQ{pde}, and $\varphi^y \leq g$ on $\Gamma_D$. It follows that $u(y) \geq \varphi^y(y) = g(y) - \ep$. Since $y\in \Gamma_D$ and $\ep > 0$ were arbitrary, it follows that $u \geq g$ on $\Gamma_D$. Moreover, since $u \geq \varphi^y$ and $u \in \USC(\bar\Omega)$, we see that $u$ is continuous at every point $y\in\Gamma_D$.

We now argue that $u$ is a subsolution of the system
\begin{equation}\label{eq:system-exist}
\left\{ \begin{aligned}
& - \Delta_\infty u - \beta |Du| = f & \mrbox{in} & \ \Omega, \\
& D_\nu u = 0 & \mrbox{on} & \ \Gamma_N. \\
\end{aligned} \right.
\end{equation}
Select a test function $\varphi\in C^2(\bar\Omega)$ and a point $x_0\in \Omega \cup \Gamma_N$ such that the map $x\mapsto (u-\varphi)(x)$ has a strict local maximum at $x=x_0$. For sufficiently small $r>0$,
\begin{equation*}
\delta(r):= (u-\varphi)(x_0) - \max_{\bar B(x_0,r)\cap \bar \Omega} (u-\varphi) > 0.
\end{equation*}
We also take $r>0$ to be small enough that $B(x_0,r) \cap \Gamma_D = \emptyset$. We may select a subsolution $w\leq u$ of the system \EQ{system-exist} for which
\begin{equation*}
u(x_0) - w(x_0) \leq \frac{1}{2}\delta(r).
\end{equation*}
It follows that the map $x\mapsto (w-\varphi)(x)$ has a local maximum at some point $y_r\in B(x_0,r)$. We deduce that
\begin{equation} \label{eq:option-1}
-\Delta^+_\infty \varphi(y_r) - \beta |D\varphi(y_r)| \leq f(y_r),
\end{equation}
or
\begin{equation} \label{eq:option-2}
y_r\in \Gamma_N \quad \mbox{and} \quad D_\nu \varphi(y_r) \leq 0.
\end{equation}
If $x_0 \in \Omega$, then for all small enough $r>0$ we must have \EQ{option-1}. In this case we may pass to the limit $r\to 0$ and use the upper semicontinuity of the map $x\mapsto \Delta_\infty^+\varphi(x)$ to deduce that
\begin{equation*}
-\Delta^+_\infty \varphi(x_0) -\beta|D\varphi(x_0)|  \leq f(x_0).
\end{equation*}
In the case that $x_0 \in \Gamma_N$, then at least one of \EQ{option-1} or \EQ{option-2} must occur for infinitely many $r>0$ along any sequence $r_j \to 0$. In this case, we may pass to limits along a subsequence to get
\begin{equation*}
-\Delta^+_\infty \varphi(x_0) -\beta|D\varphi(x_0)|  \leq f(x_0) \quad \mbox{or} \quad D_\nu \varphi(x_0) \leq 0.
\end{equation*}
We have verified that $u$ is a subsolution of \EQ{system-exist}. Since $u$ is bounded below, \LEM{local-lipschitz} implies that $u\in C(\Omega \cup \Gamma_N)$. Above we argued that $u$ is continuous at every point on $\Gamma_D$, and thus we conclude that $u\in C(\bar\Omega)$. 

We have left to show that $u$ is a supersolution of \EQ{system-exist}. Suppose that the map $x\mapsto (u - \varphi)(x)$ has a strict local minimum at a point $x_0 \in \Omega \cup \Gamma_N$, for some $\varphi \in C^2(\bar\Omega)$. Arguing indirectly, we suppose in addition that $- \Delta_\infty^- \varphi(x_0) - \beta |D \varphi(x_0)| < f(x_0)$, and that in the case $x_0\in\Gamma_N$ we have $D_\nu \varphi(x_0) < 0$. Since $\Delta_\infty^-$ is lower semicontinuous, there is a $\delta_1 > 0$ such that
\begin{equation*}
- \Delta_\infty^- \varphi - \beta |D \varphi| < f \mrbox{in} B(x_0, \delta_1) \cap \bar \Omega,
\end{equation*}
and $B(x_0, \delta_1) \cap \Gamma_D = \emptyset$. In the case $x_0\in \Gamma_N$, we may also suppose that $\delta_1> 0$ is so small that
\begin{equation*}
D_\nu \varphi <  0 \quad \mbox{in} \ B(x_0,\delta_1) \cap \Gamma_N.
\end{equation*}
In fact, we may choose a small $\delta_2 > 0$ such that
\begin{equation*}
\left\{ \begin{aligned}
& - \Delta_\infty^- \tilde \varphi - \beta |D \tilde \varphi| < f & \mbox{in} & \ B(x_0, \delta_1/2) \cap \bar \Omega,\\
& D_\nu \varphi < 0 & \mbox{on} & \ B(x_0,\delta_1/2) \cap \Gamma_N, \ \mbox{if} \ x_0 \in\Gamma_N,
\end{aligned} \right.
\end{equation*}
where $\tilde \varphi(x) := \varphi(x) - \varphi(x_0) + u(x_0) + \delta_2(\delta_1^2 / 16 -  |x - x_0|^2)$. Now define
\begin{equation*}
\tilde u(x) := \left\{ \begin{array}{ll}
\max\{ u(x), \tilde \varphi(x) \} & \mrbox{if} x \in B(x_0, \delta_1 / 2), \\
u(x) & \mrbox{otherwise.}
\end{array} \right.
\end{equation*}
Since $\tilde u$ is the maximum of two subsolutions of \EQ{system-exist} in the domain $B(x_0, \delta_1 / 2)\cap \bar\Omega$ and is equal to $u$ outside of $B(x_0, \delta_1 / 4) \cap \bar \Omega$, we see that $\tilde u$ is a subsolution of \EQ{system-exist}. Since $\tilde u(x_0) > u(x_0)$, we derive a contradiction to the definition of $u$. The proof that $u$ is a supersolution is complete.

We have constructed a maximal solution $\overline u= u \in C(\bar\Omega)$ to our boundary-value problem \EQ{pde}. To obtain the existence of a minimal solution $\underline u$, we let $w$ be the maximal solution of the problem
\begin{equation*}
\left\{ \begin{aligned}
& -\Delta_\infty w + \beta |Dw| = - f & \mbox{in} & \ \Omega,\\
& D_\nu w = 0 & \mbox{on} & \ \Gamma_N,\\
& w = -g & \mbox{on} & \ \Gamma_D,
\end{aligned} \right.
\end{equation*}
and set $\underline u = -w$. This completes the proof of \THM{existence}.
\end{proof}

\begin{remark}
Notice that Theorems \ref{thm:comparison} and \ref{thm:existence} imply the uniqueness of solutions of \EQ{pde} for \emph{generic} $f$. That is, if $g\in C(\Gamma_D)$ and $\{ f_\alpha \}_{\alpha \in \R}$ is a family functions belonging to $C(\Omega \cup \Gamma_N) \cap L^\infty(\Omega \cup \Gamma_N)$ with the property that $f_{\alpha} > f_{\beta}$ whenever $\alpha > \beta$, then there exists an at-most countable set $\Lambda \subseteq\R$ such that the boundary-value problem
\begin{equation*}
\left\{ \begin{aligned}
& -\Delta_\infty u - \beta |Dw| =  f_\alpha & \mbox{in} & \ \Omega,\\
& D_\nu u = 0 & \mbox{on} & \ \Gamma_N,\\
& u = g & \mbox{on} & \ \Gamma_D,
\end{aligned} \right.
\end{equation*}
possesses a unique solution for every $\alpha \in \R \setminus \Lambda$. For details, see the proof of \cite[Theorem 2.16]{Armstrong:preprint}. See \cite[Section 5]{Peres:2009} for a counterexample to uniqueness in the case $\beta = 0$, $\Gamma_N=\emptyset$, and $g\equiv0$.
\end{remark}

\begin{proof}[\bf Proof of \THM{stability}]
Select a test function $\varphi \in C^2(\bar \Omega)$ and a point $x_0 \in \Omega\cup \Gamma_N$ such that the map $x\mapsto (u - \varphi)(x)$ has a strict local maximum at $x \in \Omega \cup \Gamma_N$. We must show that 
\begin{equation} \label{eq:stab-1}
-\Delta^+_\infty\varphi(x_0) - \beta|D\varphi(x_0)| \leq f(x_0),
\end{equation}
or
\begin{equation} \label{eq:stab-2}
x_0 \in \Gamma_N \quad \mbox{and} \quad D_\nu \varphi(x_0) \leq 0.
\end{equation}
We proceed by showing that \EQ{stab-1} holds, under the assumption that if $x_0 \in \Gamma_N$ then \EQ{stab-2} fails. That is, in the case $x_0 \in \Gamma_N$ we assume that $D_\nu \varphi(x_0) > 0$. We may select $r> 0$ so small that either $\bar B(x_0,r) \cap \partial \Omega = \emptyset$, or $\bar B(x_0,r) \subseteq \Omega \cup \Gamma_N$ and
\begin{equation*}
D_\nu \varphi \geq 0 \quad \mbox{on} \ \bar B(x_0,r) \cap \Gamma_N.
\end{equation*}
By shrinking $r>0$ further, if necessary, we may also assume that
\begin{equation*}
(u-\varphi)(x_0) = \max_{\bar B(x_0,r/2) \cap \bar \Omega} (u -\varphi) > \max_{ \partial B(x_0,r/2) \cap \bar \Omega} (u-\varphi).
\end{equation*}
Notice that $\varphi$ is a solution of the system
\begin{equation*}
\left\{ \begin{aligned} 
& -\Delta_\infty\varphi - \beta |D\varphi| \geq h_j:= -\Delta^+_\infty \varphi - \beta_j |D\varphi| &\mbox{in} & \ B(x_0,r) \cap \Omega,\\
& D_\nu \varphi \geq 0 & \mbox{on} & \ \bar B(x_0,r) \cap \Gamma_N,
\end{aligned} \right.
\end{equation*}
and observe that the function $h_j\in \LSC(\bar B(x_0,r) \cap \bar \Omega)$.

By \THM{max-ball}, for any $0 < \ep < r/4$ we have
\begin{equation} \label{eq:stab-3}
a^-_\ep(\beta_j) S^-_\ep \varphi_\ep - a^+_\ep(\beta_j) S^+_\ep \varphi_\ep \geq \ep (h_j)_{2\ep} \quad \mbox{in} \ B(x_0,r-2\ep) \cap \bar \Omega.
\end{equation}

For sufficiently large $j$, we can find a small number $0 < \ep_j < r/3$ and a point $x_j \in B(x_0,r/2)$ such that $\ep_j \to 0$ and $x_j \to x_0$ as $j \to \infty$, and the map $x\mapsto (u_j^{\ep_j} - \varphi_{\ep_j})(x)$ attains its maximum in the ball $\bar B(x_0,r/2)$ at $x=x_j$. This implies that
\begin{multline} \label{eq:stab-4}
a^-_{\ep_j}(\beta_j) S^-_{\ep_j} \varphi_{\ep_j}(x_j) - a^+_{\ep_j}(\beta_j) S^+_{\ep_j} \varphi_{\ep_j}(x_j) \\ \leq a^-_{\ep_j}(\beta_j) S^-_{\ep_j} u_j^{\ep_j}(x_j) - a^+_{\ep_j}(\beta_j) S^+_{\ep_j} u_j^{\ep_j}(x_j).
\end{multline}
By \THM{max-ball}, we also have
\begin{equation} \label{eq:stab-5}
a^-_{\ep_j}(\beta_j) S^-_{\ep_j} u_j^{\ep_j}(x_j) - a^+_{\ep_j}(\beta_j) S^+_{\ep_j} u_j^{\ep_j}(x_j) \leq \ep_j f_j^{2\ep_j} (x_j).
\end{equation}
Combining the inequalities \EQ{stab-3}, \EQ{stab-4}, and \EQ{stab-5}, we obtain
\begin{equation*}
\left( -\Delta^+_\infty \varphi - \beta_j |D\varphi| \right)_{2\ep_j}(x_j) = (h_j)_{2\ep_j}(x_j)  \leq f^{2\ep_j}_j (x_j).
\end{equation*}
Sending $j \to \infty$ as using the lower semicontinuity of $h_j$ and the upper semicontinuity of $f$, we obtain \EQ{stab-1}, as desired.
\end{proof}

%%%%%
%%%%%
%%%%%

\bibliographystyle{plain}
\bibliography{biased}

\end{document}